\documentclass[12pt]{article}
\usepackage{amsfonts, amssymb, amsmath}
\usepackage{enumerate}
\usepackage{esint}
\newtheorem{thm}{Theorem}[section]

\newtheorem{lem}[thm]{Lemma}
\newtheorem{prop}[thm]{Proposition}

\newtheorem{rem}[thm]{Remark}

\newcommand{\f}{\frac}

\newcommand{\vc}{\infty}

\begin{document}
\title{ Weighted estimates for commutators of some singular integrals related to Schr\"odinger operators}
\author{The Anh Bui\thanks{Department of Mathematics, Macquarie University, NSW 2109, Australia and Department of Mathematics, University of Pedagogy, HoChiMinh City, Vietnam. \newline{Email: the.bui@mq.edu.au and bt\_anh80@yahoo.com}  \newline
{\it {\rm 2000} Mathematics Subject Classification:} 42B35, 35J10.
}}

\date{}

\maketitle

 \medskip

\begin{abstract}
Let $L=-\Delta +V$ with non-negative potential $V$ satisfying some appropriate reverse H\"older inequality. In this paper, we study the boundedness of the commutators of some singular integrals associated to $L$ such as Riesz transforms and fractional integrals with the new BMO functions introduced in \cite{BHS1} on the weighted spaces $L^p(w)$ where $w$ belongs to the new classes of weights introduced by \cite{BHS2}.
\end{abstract}

\section{Introduction}
Let $L=-\Delta +V$ be  the Schr\"odinger operators on $\mathbb{R}^n$ with $n\geq 3$ where the potential $V$ is in the reverse H\"older class $RH_{q}$ for some $q>n/2$, i.e., $V$ satisfies the reverse H\"older inequality
$$
\Big(\f{1}{|B|}\int_B V(y)^q dy\Big)^{1/q}\leq \f{C}{|B|}\int_B V(y)dy
$$
for all ball $B\subset \mathbb{R}^n$.\\

In this paper, we consider the following singular integrals associated to $L$:

(i) Riesz transforms $R=\nabla L^{-1/2}$ and their adjoint $R^*=L^{-1/2}\nabla$;

(ii) Fractional integrals $I_\alpha f(x)=L^{-\alpha/2}f(x)$ for $0<\alpha<n$.

In the classical case when $V=0$, it has been shown that Riesz transforms $R$ and their commutators $R_b$ with BMO functions $b$ is bounded on $L^p(w)$ for all $1<p<\vc$ and $w$ in the Muckenhoupt classes $A_p$, see for example \cite{St}. Also, the classical fractional integrals and their commutators with BMO functions $b$ are bounded from $L^p(w^p)$ to $L^q(w^q)$ for all $1<p<n/\alpha, 1/p-1/q=\alpha/n$ and $w\in A_{1+1/p'}\cap RH_q$, or equivalently $w^q\in A_{1+\f{q}{p'}}$, where $A_p$ is the Muckenhoupt class of weights, see for example \cite{MW, ST}. Recall that a non-negative and locally integrable function $w$ is said to be in the Muckenhoupt $A_p$ classes with $1\leq p<\vc$, if the following inequality holds for all balls $B\subset \mathbb{R}^n$
\begin{equation}\label{Muckenhouptweights}
\Big(\int_Bw\Big)^{1/p}\Big(\int_Bw^{-\f{1}{p-1}}\Big)^{1/p'}\leq C|B|.
\end{equation}

Recently, in \cite{BHS2}, a new class of weights  associated to Schr\"odinger operators $L$ has been introduced. According to \cite{BHS2}, the authors defined the new classes of weights $A^L_{p}=\cup_{\theta>0}A^{L,\theta}_{p} $ for $p\geq 1$, where $A^{L,\theta}_{p}$ is the set of those weights satisfying
\begin{equation}\label{classofnewweights}
\Big(\int_Bw\Big)^{1/p}\Big(\int_Bw^{-\f{1}{p-1}}\Big)^{1/p'}\leq C|B|\Big(1+\f{r}{\rho(x)}\Big)^{\theta}
\end{equation}
for all ball $B=B(x,r)$. We denote $A^L_{\vc}=\cup_{p\geq 1}A^{L}_{p}$ where the critical radius function $\rho(\cdot)$ is defined by
\begin{equation}
\rho(x)=\sup\Big\{r>0:\f{1}{r^{n-2}}\int_{B(x,r)}V\leq 1\Big\}, \ \ x\in \mathbb{R}^n,
\end{equation}
see \cite{Sh}.\\

It is easy to see that in certain circumstances the new class $A_p^L$ is larger than the Muckenhoupt class $A_p$.
The following properties hold for the new classes $A^L_p$, see \cite[Proposition 5]{BHS2}.
\begin{prop} The following statements hold:

i) $A_{p}^L\subset A_q^L$ for $1\leq p\leq q<\vc$.

ii) If $w\in A^L_p$ with $p> 1$ then there exists $\epsilon >0$ such that $w\in A_{p-\epsilon}^L$. Consequently, $A_p^L=\cup_{q<p}A_q^L$.
\end{prop}

For the ness classes $A_p^L$, the weighted norm inequalities for the some singular integrals associated to $L$ was investigated.
\begin{thm}
(a) Let $V\in RH_q$.

\ \ \ (i) If $n/2<q<n$ and $s$ is such that $1/s=1/q-1/n$, the Riesz transforms $R^*$ are bounded on $L^p(w)$ for $s'<p<\vc$ and $w\in A^L_{p/s'}$ and hence by duality $R$ is bounded on $L^p(w)$ for $1<p<s$ with $w^{-\f{1}{p-1}}\in A^L_{p'/s'}$.

\ \ \ (ii) If $q\geq n$, the Riesz transforms  $R^*$ and $R$ are bounded on $L^p(w)$ for $1<p<\vc$ and $w\in A^L_{p}$.

(b) Let $V\in RH_q$ with $q>n/2$. Then $I_\alpha$ is bounded from $L^p(w^p)$ to $L^q(w^q)$ for all $1<p<n/\alpha, 1/p-1/q=\alpha/n$ and $w^q\in A^L_{1+\f{q}{p'}}$.
\end{thm}
For the proof we refer to Theorem 3 and Theorem 4 in \cite{BHS2}.\\

Now we consider the commutators of the Riesz transforms $R$ and $R^*$ with the BMO functions $b$. It was proved in \cite{GLP} that the commutators $R_b$ and $R^*_b$ are bounded on $L^p$ here the range of $p$ depends on the potential $V$. Then the authors in \cite{BHS1} extended the classes of BMO functions to the new class $BMO_L^\theta$ with $\theta>0$ for the boundedness of the commutators $R_b$ and $R^*_b$. We would like to give a brief overview of the results in \cite{BHS1}. According to \cite{BHS1}, the new BMO space $BMO_L^\theta$ with $\theta>0$ is defined as a set of all locally integrable functions $b$ satisfying
\begin{equation}\label{eq1-intro}
\f{1}{|B|}\int_B|b(y)-b_B|dy\leq C\Big(1+\f{r}{\rho(x)}\Big)^\theta
\end{equation}
where $B=B(x,r)$ and $b_B=\f{1}{|B|}\int_B b$. A norm for $b\in BMO_L^\theta$, denoted by $\|b\|_\theta$, is given by the infimum of the constants satisfying (\ref{eq1-intro}). Clearly $BMO_L^{\theta_1}\subset BMO_L^{\theta_2}$ for $\theta_1\leq \theta_2$ and $BMO_L^0=BMO$. We define $BMO_L^\vc=\cup_{\theta>0}BMO_L^\theta$.

The following result can be considered to be a variant of John-Nirenberg inequality for the spaces $BMO_L^\theta$, see \cite[pp.119-120]{BHS1}.
\begin{prop}\label{JNforBMOL}
Let $\theta>0, s\geq 1$. If $b\in BMO_L^\theta$ then for all $B=(x_0, r)$

i) $$
\Big(\f{1}{|B|}\int_B|b(y)-b_B|^sdx\Big)^{1/s}\lesssim \|b\|_{\theta}\Big(1+\f{r}{\rho(x_0)}\Big)^{\theta '}
$$
where $\theta ' = (N_0+1)\theta$ and $N_0$ is a constant in (\ref{rho(x)and rho(y)1}).

ii) $$
\Big(\f{1}{|2^kB|}\int_{2^kB}|b(y)-b_B|dx\Big)^{1/s}\lesssim \|b\|_{\theta} k\Big(1+\f{2^kr}{\rho(x_0)}\Big)^{\theta '}
$$
for all $k\in \mathbb{N}$.
\end{prop}

Let $T$ be either $R$ or $R^*$. For $b\in BMO_L^\vc$ we consider the commutators
$$
T_bf(x)=T(bf)(x)-b(x)Tf(x)
$$
and
$$
I^b_\alpha f(x)=I_\alpha(bf)(x)-b(x)I_\alpha f(x).
$$
It was proved in \cite[Theorem 1]{BHS1} that
\begin{thm}
Let $b\in BMO_L^\theta$ with $\theta>0$ and $V\in RH_q$.

(i) If $n/2<q<n$ and $s$ is such that $1/s=1/q-1/n$, the commutators $R^*_b$ are bounded on $L^p$ for $s'<p<\vc$ and hence by duality $R_b$ is bounded on $L^p$ for $1<p<s$.

(ii) If $q\geq n$, the commutators  $R^*_b$ and $R_b$ are bounded on $L^p$ for $1<p<\vc$..
\end{thm}

The aim of this paper is investigating the boundedness of the commutators $R^*_b$, $R_b$ and $I_\alpha^b$ with $b\in BMO_L^\vc$ on the new weighted spaces $L^p(w)$ with $w\in A^L_p$.

The organization of this paper is as follows. In Section 2, we recall some basic properties of the critical radius function $\rho(\cdot)$ and consider weighted estimates for some localized operators. Section 3 is devoted to prove the main results on weighted estimates of the commutators $R^*_b$, $R_b$ and $I_\alpha^b$.

Finally, we make some conventions. Throughout the whole paper, we denote by C a positive constant
which is independent of the main parameters, but it may vary from line to line. The
symbol $X \lesssim Y$ means that there exists a positive constant $C$ such that $X \leq CY$.\\

\medskip

Recently, we have learned that the $A_p^L$ weighted norm inequalities for the commutators of the Riesz transforms was obtained independently in \cite{BHS3}. However, the approach in our paper is different from that in \cite{BHS3}. Moreover, the weighted norm inequalities for the commutators of fractional integrals $L^{-\alpha/2}$ is unique.

\section{Weighted estimates for some localized operators}

We would like to recall some important properties concerning the critical radius function which will play an important role to obtain the main results in the sequel, see \cite{Sh, DZ} respectively.

\begin{prop}\label{prop1}
If $V\in RH_{n/2}$, there exist $c_0$ and $N_0\geq 1$ such that
 \begin{equation}\label{rho(x)and rho(y)1}
 c_0^{-1}\rho(x)\Big(1+\f{|x-y|}{\rho(x)}\Big)^{-N_0}\leq \rho(y)\leq  c_0\rho(x)\Big(1+\f{|x-y|}{\rho(x)}\Big)^{\f{N_0}{N_0+1}}
 \end{equation}
for all $x,y\in \mathbb{R}^n$.
\end{prop}

A ball of the form $B(x,\rho(x))$ is called {\it a critical ball}. From the inequality (\ref{rho(x)and rho(y)1}), we can imply that for $x,y\in \sigma Q$ where $Q$ is a critical ball and $\sigma>0$, then
\begin{equation}\label{rho(x)and rho(y)2}
 \rho(x)\leq  c_\sigma\rho(y)
 \end{equation}
where $c_\sigma = c_0^2(1+\sigma)^{\f{2N_0+1}{N_0+1}}$.

\begin{prop}\label{coveringlemma}
There exists a sequence of points $x_j, j\geq 1$ in $\mathbb{R}^n$ so that the family $Q_j:=B(x_j, \rho(x_j))$ satisfies

(i) $\cup_j Q_j = \mathbb{R}^n$.

(ii) For every $\sigma \geq 1$ there exist constants $C$ and $N_1$ such that $\sum_j \chi_{\sigma Q_j}\leq C\sigma^{N_1}$.
\end{prop}

Following \cite{BHS1}, we introduce the following maximal functions for $g\in L^1_{{\rm loc}}(\mathbb{R}^n)$ and $x\in \mathbb{R}^n$
$$
M_{\rho, \alpha}g(x)=\sup_{x\in B\in \mathcal{B}_{\rho,\alpha}}\f{1}{|B|}\int_B|g|,
$$
$$
M^\sharp_{\rho, \alpha}g(x)=\sup_{x\in B\in \mathcal{B}_{\rho,\alpha}}\f{1}{|B|}\int_B|g-g_B|,
$$
where $\mathcal{B}_{\rho,\alpha}=\{B(y,r): y\in \mathbb{R}^n \ \text{and} \ r\leq \alpha\rho(y) \}$.

Also, given a ball $Q$, we define the following maximal functions for $g\in L^1_{{\rm loc}}(\mathbb{R}^n)$ and $x\in Q$
$$
M_{Q}g(x)=\sup_{x\in B\in \mathcal{F}(Q)}\f{1}{|B\cap Q|}\int_{B\cap Q}|g|,
$$
$$
M^\sharp_{Q}g(x)=\sup_{x\in B\in F(Q)}\f{1}{|B\cap Q|}\int_{B\cap Q}|g-g_{B\cap Q}|,
$$
where $\mathcal{F}(Q)=\{B(y,r): y\in Q, r>0 \}$.

We have the following lemma.
\begin{lem}\label{FSinequalityversion}
For $1<p<\vc$, then there exist $\beta$ and $\gamma$ such that if $\{Q_k\}_{k}$ is a sequence of balls as in Proposition \ref{coveringlemma} then
$$
\int_{\mathbb{R}^n}|M_{\rho,\beta}g(x)|^pw(x)dx \lesssim \int_{\mathbb{R}^n}|M^\sharp_{\rho,\gamma}g(x)|^pw(x)dx+\sum_{k}w(Q_k)\Big(\f{1}{|2Q_k|}\int_{2Q_k}|g|\Big)^p
$$
for all $g\in L^1_{{\rm loc}}(\mathbb{R}^n)$ and $w\in A^L_{\vc}.$
\end{lem}
\emph{Proof:} Note that the unweighted estimate of Lemma \ref{FSinequalityversion} was obtained \cite[Lemma 2]{BHS1}, and the weighted estimate was obtained in \cite{B} for the particular case $\rho =1$. Now we adapt some ideas in \cite[Lemma 2]{BHS1} (see also \cite{B}) to our present setting.

According to \cite[p. 121]{BHS1}, there exists $\beta>0$ so that for all critical balls $Q$ and $x\in Q$, we have
$$
M_{\rho,\beta}g(x)\leq M_{2Q}(g\chi_{2Q})(x),
$$
and for $x\in 2Q$,
$$
M^\sharp_{2Q}(g\chi_{2Q})(x)\leq M^\sharp_{\rho,2}g(x).
$$
Therefore, by the similar argument to that in \cite[Lemma 2.4]{B} we obtain for $w\in A_\vc^L$,
\begin{equation*}
\begin{aligned}
\int_{\mathbb{R}^n}|M_{\rho,\beta}g(x)|^pw(x)dx &\leq \sum_k \int_{Q_k}|M_{\rho,\beta}g(x)|^pw(x)dx\\
&\leq \sum_k \int_{Q_k}|M_{2Q}(g\chi_{2Q})(x)|^pw(x)dx\\
&\lesssim \sum_{k}\int_{2Q_k}|M^\sharp_{2Q_k}g(x)|^pw(x)dx+\sum_{k}w(Q_k)\Big(\f{1}{|2Q_k|}\int_{2Q_k}|g|\Big)^p\\
&\lesssim \sum_{k}\int_{2Q_k}|M^\sharp_{\rho,2}g(x)|^pw(x)dx+\sum_{k}w(Q_k)\Big(\f{1}{|2Q_k|}\int_{2Q_k}|g|\Big)^p\\
&\lesssim \int_{\mathbb{R}^n}|M^\sharp_{\rho,2}g(x)|^pw(x)dx+\sum_{k}w(Q_k)\Big(\f{1}{|2Q_k|}\int_{2Q_k}|g|\Big)^p.
\end{aligned}
\end{equation*}
This completes our proof.

Throughout this paper, we always assume that $N$ is a sufficiently large number and different from line to line. For $\kappa\geq 1, 0<\alpha<n$ and $1\leq s<n/\alpha$, we define the following functions for $g\in L^1_{{\rm loc}}(\mathbb{R}^n)$  and $x\in \mathbb{R}^n$
$$
G_{\kappa,\alpha,s}f(x)=\sup_{Q\ni x; Q \ {\rm is \ critical}}\sum_{k=0}^\vc 2^{-Nk}\Big(\f{1}{|2^k\widehat{Q}|^{1-\alpha s/n}}\int_{2^k\widehat{Q}}|f(z)|^sdz\Big)^{1/s}
$$
and
$$
H_{\kappa,s}f(x)=\sup_{Q\ni x; Q \ {\rm is \ critical}}\sum_{k=0}^\vc 2^{-Nk}\Big(\f{1}{|2^k\widehat{Q}|}\int_{2^k\widehat{Q}}|f(z)|^sdz\Big)^{1/s}
$$

where $\widehat{Q}=\kappa Q$.

When $\kappa=1$, we write $G_{\alpha, s}$ and $H_s$ instead of $G_{1,\alpha, s}$ and $H_{1,s}$, respectively. We are now in position to establish the weighted estimates for $G_{\kappa,\alpha,s}$ and $H_{\kappa,s}$.
\begin{prop}\label{weighted estiamtes for G}
(i) Let $w^q\in A_{1+\f{q/s}{(p/s)'}}^L$. If $p> s$ and $1/p-1/q=\alpha/n$, then we have
$$
\|G_{\kappa,\alpha,s}f\|_{L^q(w^q)}\lesssim \|f\|_{L^p(w^p)}.
$$

(ii) Let $w\in A_{p/s}^L$. If $p> s$, then we have
$$
\|H_{\kappa,s}f\|_{L^p(w)}\lesssim \|f\|_{L^p(w)}.
$$
\end{prop}
\emph{Proof:}

(i) Without of loss of generality, we can assume that $\kappa=1$. Assume that $Q=B(x_0,\rho(x_0))$. For $x\in Q$, the inequalities (\ref{rho(x)and rho(y)1}) tells us that
$$
C_0^{-1}\rho(x_0)\leq \rho(x)\leq C_0\rho(x_0).
$$
This implies $|B(x,\rho(x))|\approx |Q|$ and $Q\subset 2C_0 B(x,\rho(x))$. Therefore,
$$
G_{\alpha, s}f(x)\lesssim \sum_{k=0}^\vc 2^{-Nk}\Big(\f{1}{|2^kB(x,\rho(x))|^{1-\alpha s/n}}\int_{B_k(x,\rho(x))}|f(z)|^sdz\Big)^{1/s}
$$
where $B_k(x,\rho(x))=2C_0\times 2^kB(x,\rho(x))$.

Let $\{Q_j\}_j$ be the family of critical balls as in Proposition \ref{coveringlemma}. By (\ref{rho(x)and rho(y)1}), $C_0^{-1}\rho(x)\leq \rho(x_j)\leq C_0\rho(x)$ for all $x\in Q_j$. So, $|B(x,\rho(x))|\approx |Q_j|$ and $B_k(x,\rho(x))\subset Q_j^k$ where $Q_j^k= 2C_0\times 2^kQ_j$. For $w^q\in A_{1+\f{1}{(p/s)'}}$, using H\"older inequalities, we obtain
\begin{equation}\label{eq1-G}
\begin{aligned}
\|G_{\alpha, s}f\|_{L^q(w^q)}&\lesssim  \sum_{k=0}^\vc 2^{-Nk}\Big(\sum_j\int_{Q_j}\Big(\f{1}{|2^kB(x,\rho(x))|^{1-\alpha s/n}}\int_{B_k(x,\rho(x))}|f(z)|^sdz\Big)^{q/s}w^q(x)dx\Big)^{1/q}\\
&\lesssim  \sum_{k=0}^\vc 2^{-Nk}\Big(\sum_j\int_{Q_j}\Big(\f{1}{|2^kQ_j|^{1-\alpha s/n}}\int_{Q_j^k}|f(z)|^sdz\Big)^{q/s}w^q(x)dx\Big)^{1/q}\\
&\lesssim  \sum_{k=0}^\vc 2^{-Nk}\Big(\sum_j\f{w^q(Q_j)}{|2^kQ_j|^{q/s-\alpha q/n}}\Big(\int_{Q_j^k}|f(z)|^sdz\Big)^{q/s}\Big)^{1/q}\\
&\lesssim  \sum_{k=0}^\vc 2^{-Nk}\Big(\sum_j \f{1}{|2^kQ_j|^{q/s-\alpha q/n}}\\
&~~~~~~\times w^q(Q^k_j) \Big(\int_{Q_j^k}(w^q)^{-\f{s}{q}(\f{p}{s})'}\Big)^{\f{q/s}{(p/s)'}}\Big(\int_{Q_j^k}|f(z)|^pw^p(z)dz\Big)^{q/p}\Big)^{1/q}.
\end{aligned}
\end{equation}
Since $w^q\in A_{1+\f{q/s}{(p/s)'}}^L$, we have, by (\ref{rho(x)and rho(y)1}),
$$
w^q(Q^k_j) \Big(\int_{Q_j^k}(w^q)^{-\f{s}{q}(\f{p}{s})'}\Big)^{\f{q/s}{(p/s)'}}\leq C|Q_j^k|^{1/s-\alpha/n}2^{kN_0\theta\times (1/s-\alpha/n)}
$$
for some $\theta>0$.

This in combination with (\ref{eq1-G}) gives
\begin{equation*}
\begin{aligned}
\|G_{\alpha,s}f\|_{L^q(w^q)}&\lesssim \sum_{k}2^{-Nk}\Big(\sum_j\int_{Q_j^k}|f(z)|^pw^p(z)dz\Big)^{1/p}\\
&\lesssim  \sum_{k}2^{-Nk}\Big(\int_{\mathbb{R}^n}|f(z)|^pw^p(z)dz\Big)^{1/p}\\
&\lesssim C\|f\|_{L^p(w^p)}
\end{aligned}
\end{equation*}
where in the last inequality we used (ii) in Proposition \ref{coveringlemma}.

(ii) The proof of (ii) is similar to one of (i) and hence we omit details here.

For $0\leq \alpha<n$, let $M_\alpha$ be the fractional maximal function defined by
$$
M_\alpha f(x)=\sup_{B\ni x}\f{1}{|B|^{1-\alpha/n}}\int_B|f(y)|dy.
$$
For $s\geq 1$, we define
$$M_{\alpha,s}f=\sup_{B\ni x}\Big(\f{1}{|B|^{1-s\alpha/n}}\int_B|f(y)|^sdy\Big)^{1/s}$$.

For a family of balls $\{Q_k\}_k$ given by Proposition \ref{coveringlemma}, we define the operator $\widetilde{M}_{\alpha,s}$ as follows
\begin{equation}\label{defnoftidleM}
\widetilde{M}_{\alpha,s}f=\sum_k\chi_{Q_k}M_{\alpha,s}(f\chi_{\widetilde{Q}_k})
\end{equation}
where $\widetilde{Q}_j= 4(2C_0^2+1)\gamma Q_j$ and $\gamma$ is a constant in Proposition \ref{FSinequalityversion}.

\begin{rem}\label{rem1}
(i) For $s<p<\vc$ and $1/p-1/q=\alpha/n$, it was proved in \cite{MW} that $M_{\alpha,s}$ is bounded from $L^p(w^p)$ to $L^q(w^q)$ with $w^q\in A_{1+\f{(q/s)}{(p/s)'}}$. This together with \cite[Proposition 4]{BHS2} shows that $\widetilde{M}_{\alpha,s}$ is bounded from $L^p(w^p)$ to $L^q(w^q)$ with $w^q\in A^L_{1+\f{(q/s)}{(p/s)'}}$ here $s<p<\vc$ and $1/p-1/q=\alpha/n$.

(ii) When $\alpha =0$, we write $\widetilde{M}_{s}$ instead of $\widetilde{M}_{0,s}$. The similar argument as in (i) also shows that for $p>s$ and $w\in A^L_{p/s}$, $\widetilde{M}_{s}$ is bounded on $L^p(w)$.
\end{rem}

\section{Weighted estimates for commutators of singular integrals}
\subsection{Riesz transforms}
\subsubsection{Kernel estimates of Riesz transforms}
In the sequel, let us remind that for the number $N$, we shall mean that $N$ is a sufficiently large number and different from line to line.

Let $K$ and $K^*$ be the vector valued kernels of $R$ and $R^*$ respectively. The following propositions give some estimates on the kernels of $R$ and $R^*$, see for example \cite{Sh, GLP}.
\begin{prop}
a) If $V\in RH_{q}$ with $q> n/2$ then we have
\begin{enumerate}[(i)]
\item For every $N$ there exists a constant $C$ such that
    \begin{equation}\label{boundedK^*}
       |K^*(x,y)|\leq C\f{(1+\f{|x-z|}{\rho(x)})^{-N}}{|x-z|^{n-1}}\Big(\int_{B(z,|x-z|/4)}\f{V(u)}{|u-z|^{n-1}}+\f{1}{|x-z|}\Big).
    \end{equation}
\item For every $N$ and $0<\delta<\min\{1,2-n/q\}$ there exists a constant $C$ such that
    \begin{equation}\label{HoldercontinuityforK^*}
       |K^*(x,z)-K^*(y,z)|\leq C\f{|x-y|^\delta(1+\f{|x-z|}{\rho(x)})^{-N}}{|x-z|^{n-1+\delta}}\Big(\int_{B(z,|x-z|/4)}\f{V(u)}{|u-z|^{n-1}}+\f{1}{|x-z|}\Big)
    \end{equation}
    whenever $|x-y|<\f{2}{3}|x-z|$.
\end{enumerate}
b) If $V\in RH_{q}$ with $q\geq n$ then we have
\begin{enumerate}[(i)]
\item For every $N$ there exists a constant $C$ such that
    \begin{equation}\label{boundedK}
       |K(x,y)|\leq C\f{(1+\f{|x-y|}{\rho(x)})^{-N}}{|x-y|^{n}}.
    \end{equation}
\item For every $N$ and $0<\delta<\min\{1,1-d/q\}$ there exists a constant $C$ such that
    \begin{equation}\label{HoldercontinuityforK}
       |K^*(x,z)-K^*(y,z)|\leq C\f{|x-y|^\delta(1+\f{|x-z|}{\rho(x)})^{-N}}{|x-z|^{n+\delta}}
    \end{equation}
    whenever $|x-y|<\f{2}{3}|x-z|$.
\end{enumerate}
\end{prop}

\subsubsection{Commutators of Riesz transforms}

The main result concerning the weighted estimates for $R^*_b$ and $R_b$ is formulated by the following theorem.
\begin{thm}\label{commutatorRiezstransforms}
Let $b\in BMO_L^\theta$ with $\theta>0$ and $V\in RH_q$.

(i) If $n/2<q<n$ and $s$ is such that $1/s=1/q-1/n$, the commutator $R^*_b$ is bounded on $L^p(w)$ for $s'<p<\vc$ and $w\in A^L_{p/s'}$ and hence by duality $R_b$ is bounded on $L^p(w)$ for $1<p<s$ with $w^{-\f{1}{p-1}}\in A^L_{p'/s'}$.

(ii) If $q\geq n$, the commutators $R^*_b$ and $R_b$ are bounded on $L^p(w)$ for $1<p<\vc$ and $w\in A^L_{p}$.
\end{thm}
\emph{Proof:}
The proof of  part (ii) is completely analogous to  that of (i). Hence, we only provide the proof for (i) here and leave the second part to the interested readers.

(i) To prove (i), we exploit the strategy in \cite{BHS1}. For any $s'<p<\vc$ and $w\in A^L_{p/s'}$, we have by Lemma \ref{FSinequalityversion}
\begin{equation*}
\begin{aligned}
\|R^*_bf\|^p_{L^p(w)}&\leq \int_{\mathbb{R}^n}|M_{\rho,\beta}(R^*_bf)(x)|^pw(x)dx\\
&\lesssim \int_{\mathbb{R}^n}|M^\sharp_{\rho,\gamma}(R^*_bf)(x)|^pw(x)dx+\sum_{k}w(Q_k)\Big(\f{1}{|2Q_k|}\int_{2Q_k}|R^*_bf|\Big)^p
\end{aligned}
\end{equation*}
where $\{Q_k\}$ is a family of critical balls given in Lemma \ref{FSinequalityversion}.

{\bf 1. \ Estimate $\sum_{k}w(Q_k)\Big(\f{1}{|2Q_k|}\int_{2Q_k}|R^*_bf|\Big)^p$}

Let $s'<p_0<p$ and $Q=B(x_0,\rho(x_0))$. Then we write
$$
R^*_bf=(b-b_Q)R^*f-R^*((b-b_Q)f).
$$
So, we have
$$
\f{1}{|2Q|}\int_{2Q}|R^*_bf|dx\leq \f{1}{|2Q|}\int_{2Q}|(b-b_Q)R^*f|dx+\f{1}{|2Q|}\int_{2Q}|R^*((b-b_Q)f)|dx: =I_1+I_2.
$$
Let us estimate $I_1$ first. By H\"older inequality, we can write
\begin{equation*}
\begin{aligned}
I_{1}&\lesssim \|b\|_{\theta} \Big(\f{1}{|2Q|}\int_{2Q}|R^*f|^{p_0}\Big)^{1/p_0}\\
& \lesssim \|b\|_{\theta}\Big(\Big(\f{1}{|2Q|}\int_{2Q}|R^*f_1|^{p_0}\Big)^{1/p_0}+\Big(\f{1}{|2Q|}\int_{2Q}|R^*f_2|^{p_0}\Big)^{1/p_0}\Big)\\
&:=I_{11}+I_{12}
\end{aligned}
\end{equation*}
where $f=f_1+f_2$ with $f_1=f\chi_{4Q}$.

Due to $L^{p_0}$-boundedness of $R^*$, one has
\begin{equation*}
\begin{aligned}
I_{11}\lesssim \Big(\f{1}{|2Q|}\int_{4Q}|f|^{p_0}\Big)^{1/p_0}\lesssim \inf_{z\in Q}H_{p_0}f(z).
\end{aligned}
\end{equation*}

To estimate $I_{12}$, for $x\in 2Q$, due to (\ref{boundedK^*}), we have
\begin{equation*}
\begin{aligned}
R^*f_2(x)&\leq \int_{\mathbb{R}^n\backslash 4Q} |K^*(x,y)f(y)|dy\\
&\lesssim \int_{\mathbb{R}^n\backslash 4Q}\f{(1+\f{|x-y|}{\rho(x)})^{-N}}{|x-y|^{n}}|f(y)|dy\\
&~~~~~~+\int_{\mathbb{R}^n\backslash 4Q}\f{(1+\f{|x-y|}{\rho(x)})^{-N}}{|x-y|^{n-1}}\int_{B(y,|x-y|/4)}\f{V(u)}{|u-y|^{n-1}}|f(y)|dudy\\
&:=A_1(x)+A_2(x).
\end{aligned}
\end{equation*}
To take care $A_1$, note that $\rho(x)\approx \rho(x_0)$ and $|x-y|\approx |x_0-y|$ for all $x\in 2Q$ and $y\in \mathbb{R}^n\backslash 4Q$. So, decomposing  $\mathbb{R}^n\backslash 4Q$ into annuli $2^{k+1}\backslash 2^k Q$, we have
\begin{equation*}
\begin{aligned}
A_1(x)&\lesssim \sum_{k\geq 2}\f{2^{-kN}}{|2^kQ|}\int_{2^kQ}|f(y)|dy\\
&\lesssim \sum_{k\geq 2}2^{-kN}\Big(\f{1}{|2^kQ|}\int_{2^kQ}|f(y)|^{p_0}dy\Big)^{1/p_0}\\
&\lesssim \inf_{z\in Q}H_{p_0}f(z)
\end{aligned}
\end{equation*}
for all $x\in 2Q$.

For the term $A_2$, by decomposing $\mathbb{R}^n\backslash 4Q$ into annuli $2^{k+1}\backslash 2^k Q$, we get that
\begin{equation*}
\begin{aligned}
A_2(x)&\lesssim \sum_{k\geq 2}\f{2^{-kN}}{|2^kQ|^{1-1/n}}\int_{2^kQ}|f(y)|\int_{2^{k+2}Q}\f{V(u)}{|u-y|^{n-1}}dudy\\
&\lesssim \sum_{k\geq 2}\f{2^{-kN}}{|2^kQ|^{1-1/n}}\int_{2^kQ}|f(y)|\mathcal{I}_1(V\chi_{2^{k+2}Q})(y)dy
\end{aligned}
\end{equation*}
where $\mathcal{I}_1=(-\Delta)^{-1/2}$.

Let us remind that $\mathcal{I}_1$ is $L^{q_0}-L^{p_0'}$ boundedness with $1/p_0'=1/q_0-1/n$. This together with H\"older inequality gives
 \begin{equation*}
\begin{aligned}
A_2(x)&\lesssim \sum_{k\geq 2}\f{2^{-kN}}{|2^kQ|^{1-1/n}}\int_{\mathbb{R}^n}|(f\chi_{2^kQ})(y)|\mathcal{I}_1(V\chi_{2^{k+2}Q})(y)dy\\
&\lesssim \sum_{k\geq 2}\f{2^{-kN}}{|2^kQ|^{1-1/n}}\Big(\int_{2^kQ}|f|^{p_0}\Big)^{1/p_0}\|\mathcal{I}_1(V\chi_{2^{k+2}Q})\|_{L^{p_0'}}\\
&\lesssim \sum_{k\geq 2}\f{2^{-kN}}{|2^kQ|^{1-1/n-1/p_0}}\Big(\f{1}{|2^kQ|}\int_{2^kQ}|f|^{p_0}\Big)^{1/p_0}\|V\chi_{2^{k+2}Q}\|_{L^{q_0}}\\
\end{aligned}
\end{equation*}
Noting that if we choose $p_0$ to be close enough to $s'$ then $V\in RH_{q_0}$. This in combination with the fact that $V$ is a doubling measure gives
\begin{equation*}
\begin{aligned}
\|V\chi_{2^{k+2}Q}\|_{L^{q_0}}&\lesssim \f{1}{|2^kQ|^{1/q'_0}}\int_{2^{k+2}Q}V\\
&\lesssim \f{2^{k\kappa}}{|2^kQ|^{1/q'_0}}\int_{Q}V \ \ \text{for some $\kappa>0$}\\
&\lesssim \f{2^{k\kappa}}{|2^kQ|^{1/q'_0-1+2/n}}.
\end{aligned}
\end{equation*}
Hence,
\begin{equation*}
\begin{aligned}
A_2(x)&\lesssim \sum_{k\geq 2}\f{2^{k\kappa}}{|2^kQ|^{1/q'_0-1+2/n}}\f{2^{-kN}}{|2^kQ|^{1-1/n-1/p_0}}\Big(\f{1}{|2^kQ|}\int_{2^kQ}|f|^{p_0}\Big)^{1/p_0}\\
&\lesssim \sum_{k\geq 2}2^{-k(N-\kappa)}\Big(\f{1}{|2^kQ|}\int_{2^kQ}|f|^{p_0}\Big)^{1/p_0}\\
&\lesssim \inf_{z\in Q}H_{p_0}f(z)
\end{aligned}
\end{equation*}
for all $x\in 2Q$.

From the estimates of $A_1$ and $A_2$, we obtain $I_{12} \lesssim \inf_{z\in Q}H_{p_0}f(z)$.

The term $I_2$ can be estimated in the same line with $I_1$. Using the decomposition $f=f_1+f_2$ again, one gets that
$$
I_2\leq \Big(\f{1}{|2Q|}\int_{2Q}|R^*((b-b_Q)f_1)(y)|dy\Big)+\Big(\f{1}{|2Q|}\int_{2Q}|R^*((b-b_Q)f_2)(y)|dy\Big):=I_{21}+I_{22}.
$$
Choose $s'<r<p_0$. Using H\"older inequality and $L^r$-boundedness of $R^*$, we have
\begin{equation*}
\begin{aligned}
I_{21}&\lesssim \Big(\f{1}{|2Q|}\int_{2Q}|R^*((b-b_Q)f_1)(y)|^rdy\Big)^{1/r}\\
&\lesssim \Big(\f{1}{|2Q|}\int_{4Q}|((b-b_Q)f_1)(y)|^rdy\Big)^{1/r}\\
&\lesssim \Big(\f{1}{|2Q|}\int_{4Q}|f_1(y)|^{p_0}dy\Big)^{1/p_0}\Big(\f{1}{|2Q|}\int_{4Q}(b-b_Q)^{\nu}dy\Big)^{1/\nu}  \ \text{for some $\nu>r$}\\
& \lesssim \|b\|_\theta\inf_{z\in Q}H_{p_0}f(z).
\end{aligned}
\end{equation*}
The estimate of $I_{22}\lesssim \|b\|_\theta\inf_{z\in Q}H_{p_0}f(z)$ can be taken care similarly to ones of $I_{12}$ and $I_{21}$. So we omit the details here. To sum up, it had proved that for any critical ball $Q$, we have
\begin{equation}\label{eq1-Rb}
\f{1}{|2Q|}\int_{2Q}|R^*_bf|dx\lesssim \|b\|_\theta\inf_{z\in Q}H_{p_0}f(z).
\end{equation}

Return to the estimate of $\sum_{k}w(Q_k)\Big(\f{1}{|2Q_k|}\int_{2Q_k}|R^*_bf|\Big)^p$. Due to (\ref{eq1-Rb}), we have
\begin{equation*}
\begin{aligned}
\sum_{k}w(Q_k)\Big(\f{1}{|2Q_k|}&\int_{2Q_k}|R^*_bf|\Big)^p\\
&\lesssim \sum_{k}\|b\|^p_\theta w(Q_k)\Big(\inf_{z\in Q}H_{p_0}f(z)\Big)^p\\
&\lesssim \sum_{k}\|b\|^p_\theta\int_{Q_k}|H_{p_0}f(z)|^pw(z)dz\\
&\lesssim \|b\|^p_\theta\int_{\mathbb{R}^n}|H_{p_0}f(z)|^pw(z)dz\\
&\lesssim \|b\|^p_\theta\|f\|_{L^p(w)}^p \ \text{due to Proposition \ref{weighted estiamtes for G}}
\end{aligned}
\end{equation*}
for all $w\in A^L_{p/p_0}$. Letting $p_0\rightarrow s'$, we obtain
$$
\sum_{k}w(Q_k)\Big(\f{1}{|2Q_k|}\int_{2Q_k}|R^*_bf|\Big)^p\lesssim \|b\|^p_\theta\|f\|_{L^p(w)}^p
$$
for all $w\in A^L_{p/s'}$.

{\bf 2. \ Estimate $\int_{\mathbb{R}^n}|M^\sharp_{\rho,\gamma}(R^*_bf)(x)|^pw(x)dx$}

For any ball $B(x_0, r)$ with $r\leq \gamma\rho(x_0)$ and $x\in B$, we write
\begin{equation*}
\begin{aligned}
\f{1}{|B|}\int_B&|R^*_bf(x)-(R^*_bf)_B|dx\\
&\leq \f{2}{|B|}\int_B|(b-b_B)R^*f(x)|dx+\f{2}{|B|}\int_B|R^*((b-b_B)f_1)(x)|dx\\
& ~~~~~ +\f{1}{|B|}\int_B|R^*((b-b_B)f_2)(x)-(R^*((b-b_B)f_2))_B|dx\\
&:=E_1+E_2+E_3.
\end{aligned}
\end{equation*}
where $f=f_1+f_2$ with $f_1=f\chi_{2B}$.

Let $s'<p_0<p$, H\"older inequality and Proposition \ref{JNforBMOL} show that
\begin{equation*}
\begin{aligned}
 E_1&\lesssim \Big(\f{1}{|B|}\int_B|b-b_B|^{p_0'}\Big)^{1/p_0'}\Big(\f{1}{|B|}\int_B|R^*f|^{p_0}\Big)^{1/p_0}\\
 &\lesssim \|b\|_{\theta} \Big(\f{1}{|B|}\int_B|R^*f|^{p_0}\Big)^{1/p_0}.
\end{aligned}
\end{equation*}

For any critical ball $Q_j$ such that $x\in Q_j\cap B$. It can be verified that $B\subset \widetilde{Q}_j$. This yields that
$$
E_1 \lesssim \|b\|_{\theta}\times \inf_{y\in B}\widetilde{M}_{p_0}(R^*f)(y).
$$
For some $s'<r<p_0<p$, H\"older inequality and and Proposition \ref{JNforBMOL} again tell us that
\begin{equation*}
\begin{aligned}
 E_2&\lesssim \Big(\f{1}{|B|}\int_B|R^*((b-b_B)f_1)|^{r}\Big)^{1/r}\\
 &\lesssim \Big(\f{1}{|B|}\int_{2B}|(b-b_B)f_1|^{r}\Big)^{1/r}\\
&\lesssim \Big(\f{1}{|B|}\int_{2B}|(b-b_B)|^{\nu}\Big)^{1/\nu}\Big(\f{1}{|B|}\int_{2B}|f|^{p_0}\Big)^{1/p_0} \ \text{for some $\nu>r$}\\
&\lesssim \|b\|_{\theta}\times \inf_{y\in B}\widetilde{M}_{p_0}(f)(y).
\end{aligned}
\end{equation*}
To estimate $E_3$, we need the to show that
\begin{equation}\label{technicallem-eq}
\int_{\mathbb{R}^n\backslash 2B}|K^*(x,z)-K^*(y,z)||b(z)-b_B||f(z)|dz\lesssim \|b\|_\theta (\inf_{u\in B}H_{\gamma, p_0}f(u)+\inf_{u\in B}\widetilde{M}_{p_0}(f)(u))
\end{equation}
for all $f$ and $x,y\in B$.
If this holds, then we have
\begin{equation*}
\begin{aligned}
E_3&\lesssim \f{1}{|B|^2}\int_B\int_B\Big(\int_{R^n\backslash 2B}|K^*(u,z)-K^*(y,z)||b(z)-b_B||f(z)|dz\Big)dydu\\
&\lesssim \|b\|_\theta(H_{\gamma, p_0}f(x)+\widetilde{M}_{p_0}(f)(x)).
\end{aligned}
\end{equation*}
These three estimates of $E_1, E_2$ and $E_3$ give
$$
M^\sharp_{\rho,\gamma}(R^*_b)(x)\lesssim \|b\|_\theta(\widetilde{M}_{p_0}(R^*f)(x)+H_{\gamma, p_0}(x)+\widetilde{M}_{p_0}(f)(x)).
$$
This implies
$$
\|M^\sharp_{\rho,\gamma}(R^*_b)\|_{L^p(w)}\lesssim \|b\|_\theta(\|\widetilde{M}_{p_0}(R^*f)\|_{L^p(w)}+\|H_{\gamma, p_0}f\|_{L^p(w)}+\|\widetilde{M}_{p_0}(f)\|_{L^p(w)}.
$$
Since $\widetilde{M}_{p_0}, H_{\gamma, p_0}$ and $R^*f$ is bounded on $L^p(w)$ for all $w\in A^L_{p/p_0}$. Letting $p_0\rightarrow s'$, we obtain the desired results.\\

\emph{Proof of (\ref{technicallem-eq}):} We adapt some ideas of \cite[Lemma 6]{BHS1} to our present situation. Setting $Q=B(x_0,\gamma \rho(x_0))$, due to (\ref{HoldercontinuityforK^*}) and the fact that $\rho(x)\approx \rho(x_0)$ and $|x-z|\approx |x_0-z|$, we get
$$
\int_{\mathbb{R}^n\backslash 2B}|K^*(x,z)-K^*(y,z)||b(z)-b_B||f(z)|dz\lesssim K_1+K_2+K_3+K_4
$$
where
\begin{equation*}
\begin{aligned}
K_1&=r^\delta\int_{Q\backslash 2B}\f{|f(z)(b(z)-b_B)|}{|x_0-z|^{n+\delta}}dz,\\
K_2&=r^\delta\rho(x_0)^N\int_{Q^c}\f{|f(z)(b(z)-b_B)|}{|x_0-z|^{n+\delta+N}}dz,\\
K_3&=r^\delta\int_{Q\backslash 2B}\f{|f(z)(b(z)-b_B)|}{|x_0-z|^{n-1+\delta}}\int_{B(x_0,4|x_0-z|)}\f{V(u)}{|u-z|^{n-1}}dudz,\\
\end{aligned}
\end{equation*}
and
$$
K_4=r^\delta\rho(x_0)^N\int_{Q^c}\f{|f(z)(b(z)-b_B)|}{|x_0-z|^{n-1+\delta+N}}\int_{B(x_0,4|x_0-z|)}\f{V(u)}{|u-z|^{n-1}}dudz.
$$
Let $j_0$ be the smallest integer so that $2^{j_0}r\geq \gamma \rho(x_0)$. Splitting $Q\backslash 2B$ into annuli $2^{k+1}B\backslash 2^k B$ for $k=1, \ldots ,j_0$, we obtain
$$
K_1\lesssim \sum_{k=1}^{j_0}\f{2^{-k\delta}}{|2^kB|}\int_{2^kB}|f(z)(b(z)-b_B)|dz.
$$
By H\"older inequality and Proposition \ref{JNforBMOL},
\begin{equation*}
\begin{aligned}
K_1&\lesssim \sum_{k=1}^{j_0}k2^{-k\delta}\|b\|_\theta \Big(\f{1}{|2^kB|}\int_{2^kB}|f(z)|^{p_0}dz\Big)^{1/p_0}\\
\end{aligned}
\end{equation*}
Remind that if $x\in B\cap Q_j$ then $2^kB\subset \widetilde{Q}_j$ where $Q_j$ and $\widetilde{Q}_j$ are balls in (\ref{defnoftidleM}). Therefore,
$$
K_1\lesssim \sum_{k=1}^{j_0}k2^{-k\delta}\|b\|_\theta \inf_{u\in B}\widetilde{M}_{p_0}f(u).
$$

Splitting $Q^c$ into annuli and then applying H\"older inequality and Proposition \ref{JNforBMOL} again, we obtain
\begin{equation*}
\begin{aligned}
K_2&\lesssim \Big(\f{\rho(x_0)}{2^{j_0}r}\Big)^N\sum_{k=j_0-1}^{\vc}  \f{2^{-k\delta-(k-j_0)N}}{|2^{k-j_0+1} 2^{j_0-1}B|}\int_{2^kB}|f(z)(b(z)-b_B)|dz\\
&\lesssim \|b\|_\theta\Big(\f{\rho(x_0)}{2^{j_0}r}\Big)^N\sum_{k=j_0-1}^{\vc}  k2^{-k\delta-(k-j_0)(N-\theta ')}\Big(\f{1}{|2^{k-j_0+1} 2^{j_0-1}B|}\int_{2^kB}|f(z)|^{p_0}dz\Big)^{1/p_0}.
\end{aligned}
\end{equation*}
Since $2^{j_0}r\geq \gamma \rho(x_0)\geq 2^{j_0-1}r$, we get that
\begin{equation*}
\begin{aligned}
K_2 &\lesssim \|b\|_\theta\sum_{k=j_0-1}^{\vc}  2^{-(k-j_0)(N-\theta ')}\Big(\f{1}{|2^{k-j_0+1}Q|}\int_{2^{k-j_0+1}Q}|f(z)|^{p_0}dz\Big)^{1/p_0}\\
&\lesssim \|b\|_\theta\sum_{l=1}^{\vc}  2^{-l (N-\theta ')}\Big(\f{1}{|2^{l}Q|}\int_{2^{l}Q}|f(z)|^{p_0}dz\Big)^{1/p_0}\\
&\lesssim \|b\|_\theta \inf_{u\in B}H_{\gamma,p_0}f(u) \ \text{since $Q=\gamma B(x_0,\rho(x_0))$}.
\end{aligned}
\end{equation*}

It can be verified that
$$
K_3\lesssim \sum_{k=2}^{j_0}\f{2^{-k\delta}}{|2^kB|^{1-1/n}}\int_{2^kB}|f(z)(b(z)-b_B)|\mathcal{I}_1(V\chi_{2^{k+2}B})(z)dz
$$
Let $s'<q_0<p_0$, $\nu=\f{p_0q_0}{p_0-q_0}$ and $r$ such that $1/r=1/q_0'+1/n$ then by H\"older inequality and Proposition \ref{JNforBMOL}
\begin{equation*}
\begin{aligned}
K_3& \lesssim \|b\|_\theta\sum_{k=2}^{j_0}\f{k2^{-k\delta}}{|2^kB|^{-1/n+1/q'_0}} \Big(\f{1}{|2^kB|}\int_{2^jB}|f|^{p_0}\Big)^{1/p_0}\|\mathcal{I}_1(V\chi_{2^{k+2}B})\|_{L^{q_0'}}\\
& \lesssim \|b\|_\theta\sum_{k=2}^{j_0}\f{k2^{-k\delta}}{|2^kB|^{-1/n+1/q'_0}} \Big(\f{1}{|2^kB|}\int_{2^jB}|f|^{p_0}\Big)^{1/p_0}\|V\chi_{2^{k+2}B}\|_{L^{r}}\\
& \lesssim \|b\|_\theta\sum_{k=2}^{j_0}\f{k2^{-k\delta}}{|2^kB|^{-1/n+1/q'_0}} \|V\chi_{2^{k+2}B}\|_{L^{r}}\inf_{u\in B}\widetilde{M}_{p_0}f(u)
\end{aligned}
\end{equation*}
Noting that we can choose $q_0$ so that $V\in RH_{r}$, then we have
\begin{equation*}
\begin{aligned}
\|V\chi_{2^{k+2}B}\|_{L^{r}}&\lesssim \Big(\int_Q V(z)^rdz\Big)^{1/r}
\lesssim \f{1}{|Q|^{1-1/r}}\int_Q V(z)dz\\
&\lesssim \f{1}{|Q|^{2/n-1/r}}\lesssim \f{1}{|2^kB|^{2/n-1/r}}
\end{aligned}
\end{equation*}
for all $k=2,\ldots, j_0$.

So,
$$
K_3\lesssim \|b\|_\theta\sum_{k=2}^{j_0}k2^{-k\delta}\widetilde{M}_{p_0}f(x):=C\|b\|_\theta \inf_{u\in B}\widetilde{M}_{p_0}f(u).
$$
The similar arguments to ones used to obtain the estimates $K_2$ and $K_3$ gives
$$
K_4\lesssim \|b\|_\theta \inf_{u\in B}H_{\gamma,p_0}f(u).
$$
This completes our proof.

\subsection{Fractional integrals}
\subsubsection{Kernel estimates of fractional integrals}
Let $K_\alpha$ be the kernel of $I_\alpha$. The following results give the estimates on the kernel $K_\alpha(x,y)$.
\begin{prop}\label{kernelestimatesProp}
If $V\in RH_{q}$ with $q> n/2$ then we have
\begin{enumerate}[(i)]
\item For every $N$ there exists a constant $C$ such that
    \begin{equation}\label{boundedK2}
       |K_\alpha(x,y)|\leq C\f{(1+\f{|x-y|}{\rho(x)})^{-N}}{|x-y|^{n-\alpha}}.
    \end{equation}
\item There is a number $\delta>0$ such that for every $N$ there exists a constant $C$ such that
    \begin{equation}\label{HoldercontinuityforK2}
       |K_\alpha(x,y)-K_\alpha(x,z)|\leq C\f{|y-z|^\delta(1+\f{|x-y|}{\rho(x)})^{-N}}{|x-y|^{n+\delta-\alpha}}
    \end{equation}
    whenever $|y-z|<\f{1}{4}|x-y|$.
\end{enumerate}
\end{prop}

To prove Proposition \ref{kernelestimatesProp}, we need the following estimates of the heat kernels of $e^{-tL}$, see \cite[p.12]{DZ2}
\begin{prop}\label{heatkernelestimates}
Let $p_t(x,y)$ be the kernels associated to the semigroups $\{e^{-tL}\}_{t>0}$. If $V\in RH_{q}$ with $q> n/2$ then we have

(i) For every $N > 0$ there exists a constant $C$ such that
\begin{equation}\label{boundedp_t(x,y)}
       |p_t(x,y)|\leq \f{C}{t^{n/2}} \exp\Big(-c\f{|x-y|^2}{t}\Big) \Big(1+\f{\sqrt{t}}{\rho(x)}+ \f{\sqrt{t}}{\rho(y)}\Big)^{-N};
    \end{equation}
(ii) There is a number $\delta>0$ such that for every $N$ there exists a constant $C$ such that
    \begin{equation}\label{Holdercontinuityforp_t(x,y)}
    \begin{aligned}
       |p_t(x,y)-p_t(x,z)|&+|p_t(y,x)-p_t(z,x)|\\
       &\leq \f{C}{t^{n/2}} \Big(\f{|y-z|}{\sqrt{t}}\Big)^{\delta}\exp\Big(-c\f{|x-y|^2}{t}\Big) \Big(1+\f{\sqrt{t}}{\rho(x)}+ \f{\sqrt{t}}{\rho(y)}\Big)^{-N}
    \end{aligned}
    \end{equation}
    whenever $|y-z|<\f{1}{2}|x-y|$.
\end{prop}

\emph{Proof of Proposition \ref{kernelestimatesProp}:}
(i) We have, by (\ref{boundedp_t(x,y)}),
\begin{equation*}
    \begin{aligned}
       |K_\alpha(x,y)|&\leq \int_0^\infty |t^{\alpha/2}p_t(x,y)|\f{dt}{t}\\
       &\lesssim \int_0^{|x-y|^2}\f{t^{\alpha/2}}{t^{n/2}} \exp\Big(-c\f{|x-y|^2}{t}\Big) \Big(1+\f{\sqrt{t}}{\rho(x)}+ \f{\sqrt{t}}{\rho(y)}\Big)^{-N}\f{dt}{t}\\
       &~~~+\int_{|x-y|^2}^\vc\f{t^{\alpha/2}}{t^{n/2}} \exp\Big(-c\f{|x-y|^2}{t}\Big) \Big(1+\f{\sqrt{t}}{\rho(x)}+ \f{\sqrt{t}}{\rho(y)}\Big)^{-N}\f{dt}{t}\\
       &= I_1+I_2.
    \end{aligned}
\end{equation*}
Let us estimate $I_2$ first. Since $t>|x-y|^2$, we have, for $\epsilon >0$ so that $n>\alpha+\epsilon$,
   \begin{equation*}
    \begin{aligned}
       I_2&\lesssim  \int_{|x-y|^2}^\vc\f{t^{\alpha/2}}{t^{n/2}} \exp\Big(-c\f{|x-y|^2}{t}\Big) \Big(1+\f{|x-y|}{\rho(x)}\Big)^{-N}\f{dt}{t}\\
          &\lesssim  \int_{|x-y|^2}^\vc\f{t^{\alpha/2}}{t^{n/2}} (\f{t}{|x-y|^2}\Big)^{n/2-\alpha/2-\epsilon} \Big(1+\f{|x-y|}{\rho(x)}\Big)^{-N}\f{dt}{t}\\
          &\lesssim  \f{1}{|x-y|^{n-\alpha}}\Big(1+\f{|x-y|}{\rho(x)}\Big)^{-N}.
    \end{aligned}
\end{equation*}
For $I_1$, we have
\begin{equation*}
    \begin{aligned}
       I_1&\lesssim  \int^{|x-y|^2}_0\f{t^{\alpha/2}}{t^{n/2}}\Big(\f{t}{|x-y|^2}\Big)^{n/2+N/2} \Big(1+\f{\sqrt{t}}{\rho(x)}\Big)^{-N}\f{dt}{t}\\
          &\lesssim  \int^{|x-y|^2}_0\f{t^{\alpha/2}}{|x-y|^{n}}\Big(\f{\sqrt{t}}{|x-y|}\Big)^{N} \Big(\f{\sqrt{t}+\rho(x)}{\rho(x)}\Big)^{-N}\f{dt}{t}\\
          &\lesssim  \int^{|x-y|^2}_0\f{t^{\alpha/2}}{|x-y|^{n}}\Big(\f{\sqrt{t}+\rho(x)}{|x-y|+\rho(x)}\Big)^{N} \Big(\f{\sqrt{t}+\rho(x)}{\rho(x)}\Big)^{-N}\f{dt}{t}\\
          &\lesssim  \f{1}{|x-y|^{n-\alpha}}\Big(1+\f{|x-y|}{\rho(x)}\Big)^{-N}.
    \end{aligned}
\end{equation*}
This completes (i).

(ii) For $|y-z|<\f{1}{4}|x-y|$, using (\ref{Holdercontinuityforp_t(x,y)}) gives
\begin{equation*}
    \begin{aligned}
     |K_\alpha(x,y)-K_\alpha(x,z)|&\lesssim \int_0^\vc \f{t^{\alpha/2}}{t^{n/2}} \Big(\f{|y-z|}{\sqrt{t}}\Big)^{\delta}\exp\Big(-c\f{|x-y|^2}{t}\Big) \Big(1+\f{\sqrt{t}}{\rho(x)}+ \f{\sqrt{t}}{\rho(y)}\Big)^{-N}\f{dt}{t}\\
     &\lesssim \Big(\f{|y-z|}{|x-y|}\Big)^{\delta}\int_0^\vc \f{t^{\alpha/2}}{t^{n/2}}\exp\Big(-c'\f{|x-y|^2}{t}\Big) \Big(1+\f{\sqrt{t}}{\rho(x)}+ \f{\sqrt{t}}{\rho(y)}\Big)^{-N}\f{dt}{t}\\
     &\lesssim \Big(\f{|y-z|}{|x-y|}\Big)^{\delta}\Big(\int_0^{|x-y|^2}\ldots +\int_{|x-y|^2}^\vc\ldots\Big)\\
     &=\Big(\f{|y-z|}{|x-y|}\Big)^{\delta}(II_1+II_2).
    \end{aligned}
\end{equation*}
At this stage, repeating the arguments in (i), we obtain (ii).

\subsubsection{Commutators of fractional integrals}

We are now in position to state the result concerning the weighted estimates for $I_\alpha^b$.
\begin{thm}\label{fractionalIntegrals}
Let $b\in BMO_L^\theta$ with $\theta>0$ and $V\in RH_q$ with $q>n/2$. Then the commutator $I_\alpha^b$ is bounded from $L^p(w^p)$ to $L^q(w^q)$ for $1<p<\vc, 1/p-1/q=\alpha/n$  and $w^q\in A^L_{1+\f{q}{p'}}$.
\end{thm}

\emph{Proof:} The strategy of the proof for Theorem \ref{fractionalIntegrals} is similar to that of Theorem \ref{commutatorRiezstransforms}.

For any $1<s<p<\vc, 1/p-1/q=\alpha/n$ and $w^q\in A^L_{1+\f{(q/s)}{(p/s)'}}$, we have by Lemma \ref{FSinequalityversion}
\begin{equation*}
\begin{aligned}
\|I_\alpha^b f\|^q_{L^q(w^q)}&\leq \int_{\mathbb{R}^n}|M_{\rho,\beta}(I_\alpha^bf)(x)|^qw^q(x)dx\\
&\lesssim \int_{\mathbb{R}^n}|M^\sharp_{\rho,\gamma}(I_\alpha^bf)(x)|^qw^q(x)dx+\sum_{k}w^q(Q_k)\Big(\f{1}{|2Q_k|}\int_{2Q_k}|I_\alpha^bf|\Big)^q
\end{aligned}
\end{equation*}
where $\{Q_k\}$ is a family of critical balls given in Proposition \ref{FSinequalityversion}.

So, to obtain the weighted estimates for $I_\alpha^b$, we need only to consider $\int_{\mathbb{R}^n}|M^\sharp_{\rho,\gamma}(I_\alpha^bf)(x)|^qw^q(x)dx$ and $\sum_{k}w^q(Q_k)\Big(\f{1}{|2Q_k|}\int_{2Q_k}|I_\alpha^bf|\Big)^q$.

{\bf Step 1. \ Estimate $\sum_{k}w^q(Q_k)\Big(\f{1}{|2Q_k|}\int_{2Q_k}|I_\alpha^bf|\Big)^q$}

Let $1<s<p$, $1/s-1/\upsilon=\alpha/n$ and $Q=B(x_0,\rho(x_0))$. We have
$$
I_\alpha^bf=(b-b_Q)I_\alpha f-I_\alpha((b-b_Q)f).
$$
Hence,
$$
\f{1}{|2Q|}\int_{2Q}|I_\alpha^bf|dx\leq \f{1}{|2Q|}\int_{2Q}|(b-b_Q)I_\alpha f|dx+\f{1}{|2Q|}\int_{2Q}|I_\alpha((b-b_Q)f)|dx: =I_1+I_2.
$$
To take care $I_1$, using H\"older inequality and Proposition \ref{JNforBMOL}, we get that
\begin{equation*}
\begin{aligned}
I_{1}&\lesssim \|b\|_{\theta} \Big(\f{1}{|2Q|}\int_{2Q}|I_\alpha f|^{\upsilon}\Big)^{1/\upsilon}\\
& \lesssim \|b\|_{\theta}\Big(\Big(\f{1}{|2Q|}\int_{2Q}|I_\alpha f_1|^{\upsilon}\Big)^{1/\upsilon}+\Big(\f{1}{|2Q|}\int_{2Q}|I_\alpha f_2|^{\upsilon}\Big)^{1/\upsilon}\Big)\\
&:=I_{11}+I_{12}
\end{aligned}
\end{equation*}
where $f=f_1+f_2$ with $f_1=f\chi_{4Q}$.

Since $I_\alpha$  is $L^{s}-L^\upsilon$ bounded, one has
\begin{equation*}
\begin{aligned}
I_{11}\lesssim \Big(\f{1}{|2Q|^{1-\alpha s/n}}\int_{4Q}|f|^{s}\Big)^{1/s}\lesssim \inf_{z\in Q}G_{\alpha,s}f(z).
\end{aligned}
\end{equation*}

To estimate $I_{12}$, for $x\in 2Q$, (\ref{boundedK2}) implies that
\begin{equation*}
\begin{aligned}
I_\alpha f_2(x)&\leq \int_{\mathbb{R}^n\backslash 4Q} |K_\alpha(x,y)f(y)|dy\\
&\lesssim \int_{\mathbb{R}^n\backslash 4Q}\f{(1+\f{|x-y|}{\rho(x)})^{-N}}{|x-y|^{n-\alpha}}|f(y)|dy.
\end{aligned}
\end{equation*}
In this situation, we have $\rho(x)\approx \rho(x_0)$ and $|x-y|\approx |x_0-y|$ for all $x\in 2Q$ and $y\in \mathbb{R}^n\backslash 4Q$. So, decomposing  $\mathbb{R}^n\backslash 4Q$ into annuli $2^{k+1}\backslash 2^k Q$, we have, by H\"older inequality,
\begin{equation*}
\begin{aligned}
I_\alpha f_2(x)&\lesssim \sum_{k\geq 2}\f{2^{-kN}}{|2^kQ|^{1-\alpha/n}}\int_{2^kQ}|f(y)|dy\\
&\lesssim \sum_{k\geq 2}2^{-kN}\Big(\f{1}{|2^kQ|^{1-\alpha s/n}}\int_{2^kQ}|f(y)|^{s}dy\Big)^{1/s}\\
&\lesssim \inf_{z\in Q}G_{\alpha, s}f(z)
\end{aligned}
\end{equation*}
for all $x\in 2Q$.
Hence $I_{12} \lesssim \inf_{z\in Q}G_{\alpha, s}f(z)$.

The estimate for  $I_2$ can be proceeded in the same line with one of $I_1$. The decomposition $f=f_1+f_2$ gives
$$
I_2\leq \Big(\f{1}{|2Q|}\int_{2Q}|I_\alpha((b-b_Q)f_1)(y)|dy\Big)+\Big(\f{1}{|2Q|}\int_{2Q}|I_\alpha((b-b_Q)f_2)(y)|dy\Big):=I_{21}+I_{22}.
$$
Choose $1 <r<s<p$ and $1/r-1/r_0=\alpha/n$. Using H\"older inequality, Proposition \ref{JNforBMOL} and $L^r-L^{r_0}$-boundedness of $I_\alpha$, we have
\begin{equation*}
\begin{aligned}
I_{21}&\lesssim \Big(\f{1}{|2Q|}\int_{2Q}|I-\alpha((b-b_Q)f_1)(y)|^{r_0}dy\Big)^{1/r_0}\\
&\lesssim \f{1}{|2Q|^{-\alpha/n}}\Big(\f{1}{|2Q|}\int_{4Q}|((b-b_Q)f_1)(y)|^rdy\Big)^{1/r}\\
&\lesssim \f{1}{|2Q|^{-\alpha/n}}\Big(\f{1}{|2Q|}\int_{4Q}|f_1(y)|^{s}dy\Big)^{1/s}\Big(\f{1}{|2Q|}\int_{4Q}(b-b_Q)^{\nu}dy\Big)^{1/\nu}  \ \text{for some $\nu>r$}\\
& \lesssim \|b\|_\theta\inf_{z\in Q}G_{\alpha,s}f(z).
\end{aligned}
\end{equation*}
The estimate $I_{22}\lesssim \|b\|_\theta\inf_{z\in Q}G_{\alpha,s}f(z)$ can be obtained by the similar approach to ones of $I_{12}$ and $I_{21}$. So we omit the details here.

To sum up, for any critical ball $Q$, we have
\begin{equation}\label{eq1-Rb}
\f{1}{|2Q|}\int_{2Q}|I_\alpha^bf|dx\lesssim \|b\|_\theta\inf_{z\in Q}G_{\alpha,s}f(z).
\end{equation}

Return to the estimate of $\sum_{k}w^q(Q_k)\Big(\f{1}{|2Q_k|}\int_{2Q_k}|I_\alpha^b|\Big)^q$. Due to (\ref{eq1-Rb}), we have
\begin{equation*}
\begin{aligned}
\sum_{k}w^q(Q_k)\Big(\f{1}{|2Q_k|}&\int_{2Q_k}|I_\alpha^bf|\Big)^q\\
&\lesssim \sum_{k}\|b\|^q_\theta w^q(Q_k)\Big(\inf_{z\in Q}G_{\alpha, s}f(z)\Big)^q\\
&\lesssim \sum_{k}\|b\|^q_\theta \int_{Q_k}|G_{\alpha, s}f(z)|^qw^q(z)dz\\
&\lesssim \|b\|^q_\theta\int_{\mathbb{R}^n}|G_{\alpha, s}f(z)|^qw^q(z)dz\\
&\lesssim \|b\|^q_\theta\|f\|_{L^p(w^p)}^q \ \text{due to Proposition \ref{weighted estiamtes for G}}
\end{aligned}
\end{equation*}
for all $w^q\in A^L_{1+\f{q/s}{(p/s)'}}$. Letting $s\rightarrow 1$, we obtain
$$
\sum_{k}w^q(Q_k)\Big(\f{1}{|2Q_k|}\int_{2Q_k}|I_\alpha^b|\Big)^q\lesssim \|b\|^p_\theta\|f\|_{L^p(w^p)}^q
$$
for all $w^q\in A^L_{1+\f{q}{p'}}$.

{\bf Step 2. \ Estimate $\int_{\mathbb{R}^n}|M^\sharp_{\rho,\gamma}(I_\alpha^bf)(x)|^qw^q(x)dx$}

For any ball $B(x_0, r)$ with $r\leq \gamma\rho(x_0)$ and $x\in B$, we write
\begin{equation*}
\begin{aligned}
\f{1}{|B|}\int_B&|I_\alpha^bf(x)-(I_\alpha^b f)_B|dx\\
&\leq \f{2}{|B|}\int_B|(b-b_B)I_\alpha f(x)|dx+\f{2}{|B|}\int_B|I_\alpha((b-b_B)f_1)(x)|dx\\
& ~~~~~ +\f{1}{|B|}\int_B|I_\alpha((b-b_B)f_2)(x)-(I_\alpha((b-b_B)f_2))_B|dx\\
&:=E_1+E_2+E_3.
\end{aligned}
\end{equation*}
where $f=f_1+f_2$ with $f_1=f\chi_{2B}$.

Applying H\"older inequality and Proposition \ref{JNforBMOL}, we have
\begin{equation*}
\begin{aligned}
 E_1&\lesssim \Big(\f{1}{|B|}\int_B|b-b_B|^{s'}\Big)^{1/s'}\Big(\f{1}{|B|}\int_B|I_\alpha f|^{s}\Big)^{1/s
 }\\
 &\lesssim \|b\|_{\theta} \Big(\f{1}{|B|}\int_B|I_\alpha f|^{s}\Big)^{1/s}.
\end{aligned}
\end{equation*}

For any critical ball $Q_j$ such that $x\in Q_j\cap B$, it is easy to see that $B\subset \widetilde{Q}_j$. Therefore,
$$
E_1 \lesssim \|b\|_{\theta}\times \inf_{u\in B}\widetilde{M}_{s}(I_\alpha f)(u).
$$
For some $1<r<s<p$ and $1/r-1/r_0=\alpha/n$, H\"older inequality and Proposition \ref{JNforBMOL} again tell us that
\begin{equation*}
\begin{aligned}
 E_2&\lesssim \Big(\f{1}{|B|}\int_B|I_\alpha((b-b_B)f_1)|^{r_0}\Big)^{1/r_0}\\
 &\lesssim \f{1}{|B|^{-\alpha/n}}\Big(\f{1}{|B|}\int_{2B}|(b-b_B)f_1|^{r}\Big)^{1/r}\\
&\lesssim \f{1}{|B|^{-\alpha/n}}\Big(\f{1}{|B|}\int_{2B}|(b-b_B)|^{\nu}\Big)^{1/\nu}\Big(\f{1}{|B|}\int_{2B}|f|^{s}\Big)^{1/s} \ \text{for some $\nu>r$}\\
&\lesssim \|b\|_{\theta}\times \inf_{u\in B}\widetilde{M}_{\alpha, s}(f)(u).
\end{aligned}
\end{equation*}

Before taking care $E_3$, we need the to show that
\begin{equation}\label{technicallem2-eq}
\int_{\mathbb{R}^n\backslash 2B}|K_\alpha(x,z)-K_\alpha(y,z)||b(z)-b_B||f(z)|dz\lesssim \|b\|_\theta (\inf_{u\in B}G_{\gamma,\alpha, s}f(u)+\inf_{u\in B}\widetilde{M}_{\alpha,s}(f)(u))
\end{equation}
for all $f$ and $x,y\in B$.

If this holds, then we have
\begin{equation*}
\begin{aligned}
E_3&\lesssim \f{1}{|B|^2}\int_B\int_B\Big(\int_{R^n\backslash 2B}|K_\alpha(u,z)-K_\alpha(y,z)||b(z)-b_B||f(z)|dz\Big)dydu\\
&\lesssim \|b\|_\theta(G_{\gamma, \alpha, s}(x)+\widetilde{M}_{\alpha, s}(f)(x)).
\end{aligned}
\end{equation*}
These three estimates of $E_1, E_2$ and $E_3$ give that
$$
M^\sharp_{\rho,\gamma}(I_\alpha^b)(x)\lesssim \|b\|_\theta(\widetilde{M}_{\alpha,s}(f)(x)+G_{\gamma,\alpha,s}f(x)+\widetilde{M}_{s}(I_\alpha f)(x)).
$$
This implies
$$
\|M^\sharp_{\rho,\gamma}(I_\alpha^b)\|_{L^q(w^q)}\lesssim \|b\|_\theta(\|\widetilde{M}_{\alpha, s}(f)\|_{L^q(w^q)}+\|G_{\gamma,\alpha, s}f\|_{L^q(w^q)}+\|\widetilde{M}_{s}(I_\alpha f)\|_{L^q(w^q)}).
$$
Since $\widetilde{M}_{\alpha, s}$ and $G_{\gamma, \alpha, s}$ are bounded form $L^p(w^p)$ to $L^q(w^q)$ for all $1<p<n/\alpha$, $1/p-1/q=\alpha/n$ and $w^q\in A^L_{1+(q/s)/(p/s)'}$, we have
$$
\|\widetilde{M}_{\alpha, s}(f)\|_{L^q(w^q)}+\|G_{\gamma,\alpha, s}f\|_{L^q(w^q)}\lesssim \|f\|_{L^p(w^p)}.
$$
For the last term $\|\widetilde{M}_{s}(I_\alpha f)\|_{L^q(w^q)}$, from the weighted estimates of $I_\alpha$ and $\widetilde{M}_s$ (see Remark \ref{rem1}) and the fact that $A_{1+\f{(q/s)}{(p/s)'}}\subset A_{q/s}$, one gets that
$$
\|\widetilde{M}_{s}(I_\alpha f)\|_{L^q(w^q)}\lesssim \|I_\alpha f\|_{L^q(w^q)}\lesssim \|f\|_{L^p(w^p)}.
$$
Hence,
$$
\|M^\sharp_{\rho,\gamma}(I_\alpha^b)\|_{L^q(w^q)}\lesssim \|b\|_\theta\|f\|_{L^p(w^p)}
$$
for all $w^q\in A_{1+\f{(q/s)}{(p/s)'}}$ and $1<s<p$.

Letting $s\rightarrow 1$, we obtain the desired results.\\

\emph{Proof of (\ref{technicallem2-eq}):} Setting $Q=B(x_0,\gamma \rho(x_0))$, due to (\ref{HoldercontinuityforK2}) and the fact that $\rho(x)\approx \rho(x_0)$ and $|x-z|\approx |x_0-z|$, we get
$$
\int_{\mathbb{R}^n\backslash 2B}|K_\alpha(x,z)-K_\alpha(y,z)||b(z)-b_B||f(z)|dz\lesssim K_1+K_2
$$
where
$$
K_1=r^\delta\int_{Q\backslash 2B}\f{|f(z)(b(z)-b_B)|}{|x_0-z|^{n+\delta-\alpha}}dz,
$$
and
$$
K_2=r^\delta\rho(x_0)^N\int_{Q^c}\f{|f(z)(b(z)-b_B)|}{|x_0-z|^{n+\delta+N-\alpha}}dz.
$$

Let $j_0$ be the smallest integer so that $2^{j_0}r\geq \gamma \rho(x_0)$. Splitting $Q\backslash 2B$ into annuli $2^{k+1}B\backslash 2^k B$ for $k=1, \ldots ,j_0$, we obtain
$$
K_1\lesssim \sum_{k=1}^{j_0}\f{2^{-k\delta}}{|2^kB|^{1-\alpha/n}}\int_{2^kB}|f(z)(b(z)-b_B)|dz.
$$
Using H\"older inequality and Proposition \ref{JNforBMOL}, we obtain
\begin{equation*}
\begin{aligned}
K_1&\lesssim \sum_{k=1}^{j_0}k2^{-k\delta}\|b\|_\theta \Big(\f{1}{|2^kB|^{1-\alpha s/n}}\int_{2^kB}|f(z)|^{s}dz\Big)^{1/s}\\
\end{aligned}
\end{equation*}
Note that if $x\in B\cap Q_j$ then $2^kB\subset \widetilde{Q}_j$ where $Q_j$ and $\widetilde{Q}_j$ are balls in (\ref{defnoftidleM}). Therefore,
$$
K_1\lesssim \sum_{k=1}^{j_0}k2^{-k\delta}\|b\|_\theta \inf_{u\in B}\widetilde{M}_{\alpha, s}f(u).
$$

Splitting $Q^c$ into annuli and then applying H\"older inequality and Proposition \ref{JNforBMOL} again, we obtain
\begin{equation*}
\begin{aligned}
K_2&\lesssim \Big(\f{\rho(x_0)}{2^{j_0}r}\Big)^N\sum_{k=j_0-1}^{\vc}  \f{2^{-k\delta-(k-j_0)N}}{|2^{k-j_0+1} 2^{j_0-1}B|^{1-\alpha/n}}\int_{2^kB}|f(z)(b(z)-b_B)|dz\\
&\lesssim \|b\|_\theta\Big(\f{\rho(x_0)}{2^{j_0}r}\Big)^N\sum_{k=j_0-1}^{\vc}  k2^{-k\delta-(k-j_0)(N-\theta ')}\Big(\f{1}{|2^{k-j_0+1} 2^{j_0-1}B|^{1-\alpha s/n}}\int_{2^kB}|f(z)|^{s}dz\Big)^{1/s}.
\end{aligned}
\end{equation*}
Since $2^{j_0}r\geq \gamma \rho(x_0)\geq 2^{j_0-1}r$, we get that
\begin{equation*}
\begin{aligned}
K_2 &\lesssim \|b\|_\theta\sum_{k=j_0-1}^{\vc}  2^{-(k-j_0)(N-\theta ')}\Big(\f{1}{|2^{k-j_0+1}Q|^{1-\alpha s/n}}\int_{2^{k-j_0+1}Q}|f(z)|^{s}dz\Big)^{1/s}\\
&\lesssim \|b\|_\theta\sum_{l=1}^{\vc}  2^{-l (N-\theta ')}\Big(\f{1}{|2^{l}Q|^{1-\alpha s/n}}\int_{2^{l}Q}|f(z)|^{s}dz\Big)^{1/s}\\
&\lesssim \|b\|_\theta \inf_{u\in B}G_{\gamma,\alpha, s}f(u) \ \text{since $Q=\gamma B(x_0,\rho(x_0))$}.
\end{aligned}
\end{equation*}

This completes our proof.

\emph{Acknowledgements} \  \ The author would like to thank his supervisor, Prof. X. T. Duong for helpful
comments and suggestions.


\begin{thebibliography}{99999999}
\bibitem[BD]{BD} T. A. Bui and X. T. Duong, On $A_p^L$ weighted estimates for spectral multipliers of Schr\"odinger operators, preprint.

\bibitem[B]{B} T. A. Bui Weighted norm inequalities for pseudo-differential operators and
their commutators, preprint.

\bibitem[BHS1]{BHS1} B. Bongioanni, E. Haboure and O. Salinas, Commutators of Riezs transforms related to Schr\"odinger operators, J. Fourier Anal. Appl. 17 (2011), 115-134.

\bibitem[BHS2]{BHS2} B. Bongioanni, E. Haboure and O. Salinas, Classes of weights related to Schr\"odinger operators, J. Math. Anal. Appl. 373 (2011), 563-579.

\bibitem[BHS3]{BHS3} B. Bongioanni, E. Haboure and O. Salinas, Weighted inequalities for commutators of Schr\"odinger-Riesz transforms, preprint http://arxiv.org/abs/1110.5797.

\bibitem[DZ1]{DZ} J. Dziuba\'nski and J. Zienkiewicz, Hardy spaces $H^1$ associated to Schr\"odinger operators with potential satisfying reverse H\"older inequality, Rev. Mat. Iberoam. 15 (1999), 279-296.

\bibitem[DZ2]{DZ2} J. Dziuba\'nski and J. Zienkiewicz, Hardy spaces $H^p$ associated to Schr\"odinger operators with potential from reverse H\"older inequality, Colloq. Math. 98 (2003), 5-38.

\bibitem[GLP]{GLP} Z. Guo, P. Li and L. Peng, $L^p$ boundedness of commutators of Riesz transforms associated to Schr\"odinger operators, J. Math. Anal. Appl. 341 (2008), 421-432.

\bibitem[Sh]{Sh} Z. Shen, $L^p$ estimates for Schr\"odinger operators with certain potentials, Ann. Inst. Fourier (Grenoble) 45 (1995), 513-546.

\bibitem[St]{St} E.M. Stein,  {\it Harmonic analysis: Real variable
methods, orthogonality and oscillatory integrals}, Princeton Univ.
Press, Princeton, NJ, (1993).

\bibitem[MW]{MW} B. Muckenhoupt and R. Wheeden, Weighted norm inequalities for fractional integrals, Trans.
Amer. Math. Soc. 192 (1974), 261-274.

\bibitem[ST]{ST} C. Segovia and J. L. Torrea, Weighted inequalities for commutators
of fractional and singular integrals, Publ. Mat. 35 (1991), 209–235.

\end{thebibliography}
\end{document}